\documentclass[10pt]{article}
\usepackage{latexsym,epsf,amssymb,times}

\newlength{\pwdth}

\newcommand{\frpar}[1]{\begin{center}\fbox{\begin{minipage}[c]{\pwdth}#1
\end{minipage}}\end{center}}

\newenvironment{proof}{\begin{trivlist}
                       \item[]\hspace{0.3cm}{\bf Proof}
                       \hspace{0cm} }{\hfill $\Box$
                       \end{trivlist}}

\title{Spectrum through pseudospectrum}

\author{Ioannis Koutis}
\begin{document}

\maketitle

\begin{abstract}

\noindent {\bf This report was written in 2001 and it is a translation of work that was
originally published in Greek, in the author's diploma thesis in July 1998. It may
contain minor mistakes and should not be considered a complete study. It however touches
upon several of the considerations that will be included in the complete paper. The work
was done jointly with E. Gallopoulos and was presented at the FOCM99 conference, at the
5th IMACS conference on iterative methods in scientific computing, and at the 50th annual
meeting
of SIAM.} \\

\noindent We present an iteration for the computation of simple eigenvalues using a
pseudospectrum approach. The most appealing characteristic of the proposed iteration is
that it reduces the computation of a single eigenvalue to a small number of eigenvalue
computations on Hermitian matrices. We show that this number is directly associated with
the matrix pseudospectrum. We present numerical results and we discuss advantages and
drawbacks of the method. We also discuss its relationship with an iteration that was
proposed independently in [Stewart, O' Leary, ETNA, Vol 8, 1998].

\noindent

\end{abstract}

\section{A pseudospectrum setting}\label{sec:SETTING}

\newtheorem{theorem}{Theorem}

We shortly describe some well known characteristics of matrix
pseudospectrum.  Denote with $ \Lambda(A) $ the spectrum of a matrix A, with
$\Lambda_{\epsilon}(A)$ the $\epsilon$-pseudospectrum of A, with
$\sigma_{\min}=\sigma_{N}\leq \sigma_{N-1} \ldots \leq \sigma_{1}$ the
singular values of A, and with $D(z,\varrho)$ and $D^{\circ}(z,\varrho)$ a
closed disk and an open disk respectively, with center $z$ and radius
$\varrho$. Let $N$ denote a normal matrix.

\begin{theorem} \label{th:PS.BASICS} \qquad \qquad \qquad
\begin{itemize}
\item $\Lambda(A) \subset \Lambda_{\epsilon}(A)$ for $\epsilon>0$.
\item $\epsilon<\epsilon_{1} \Leftrightarrow
\Lambda_{\epsilon}(A) \subset \Lambda_{\epsilon_{1}}(A)$.
\item $\Lambda(N) = \Lambda(A)$ $\Rightarrow$  $\Lambda_{\epsilon}(N)\subseteq
\Lambda_{\epsilon}(A)$.
\end{itemize}
\end{theorem}

These properties state that $\epsilon$-pseudospectrum forms closed
curves around eigenvalues. The crucial property of subharmonicity
of the norm of the resolvent $||(zI-A)^{-1}||$ (see
\cite{Gallestey.98}), assures that
\begin{eqnarray}
\sigma_{\min}(zI-A) = \textnormal{local minimum} =0 \Leftrightarrow
z\in\Lambda(A)
\end{eqnarray}
This alternative characterization of an eigenvalue together with
the properties of matrix pseudospectrum appeal for the
corresponding optimization problem:
\begin{eqnarray*}
\min f(z)=\sigma_{\min}(zI-A), ~~z\in \mathbb{C}
\end{eqnarray*}

Notice that in this formulation we are seeking for a minimizer \textit{over
the complex plane}. This is not usually the case in other methods for the
computation of eigenvalues . The "heart" of any competent minimization
algorithm can naturally be the following theorem by Sun \cite{Sun.88}.

\begin{theorem} \label{th:gradient}
Denote with $\sigma_i$ the $i^{th}$ singular value and $u_{i}$, $v_{i}$ the
corresponding left and right singular vectors. Let also $z=x+iy \in
\mathbb{C}\setminus\Lambda(A)$. Then $g(x,y)=\sigma_i(zI-A)$ is real
analytic in a neighborhood of $(x,y)$, if $\sigma_{\min}(zI-A)$ is a simple
singular value. For the gradient of $g(x,y)$ we have
\begin{eqnarray*}
\nabla g(x,y) = (Re(v_{i}^{*}u_{i}),Im(v_{i}^{*}u_{i}))
\end{eqnarray*}
\end{theorem}

What is interesting, is that $\nabla \sigma_{\min}(zI-A)$ can be computed
almost as easy as $\sigma_{\min}$ itself. Given that $\sigma_{\min}$ can be
computed by means of some of the well known methods for Hermitian matrices,
we can readily setup the framework for a minimization algorithm. Using for
example a steepest descent strategy we get the extremely compact
{\tt{Matlab}} routine shown below

\setlength{\pwdth}{4.5in} \frpar{
\begin{tabular} [h]{l}
{\bf  Eigenvalue via Hermitian Computations}
\\
\\  {\tt function z=sdeig(A,starting\_point)}
\\
\\  {\tt z= starting\_point; s\_min= 1;}
\\  {\tt while (s\_min $>$ tol)}
\\  {\tt \qquad [u\_{min},s\_{min},v\_{min}] = svds(z*I-A);}
\\  {\tt \qquad direction = v\_{min}$^{'}*$u\_{min} /
$|$v\_{min}$^{'}*$u\_{min}$|$;}
\\  {\tt \qquad step\_size = ? ;}
\\  {\tt \qquad z = z - step\_size*direction;}
\\ {\tt end}
\end{tabular}
}

The determination of the step size is obviously of critical importance for
the effectiveness of the {\tt sdeig}. In our experiments we used the general
purpose minimization function {\tt fminu} of {\tt Matlab Optimization
Toolbox}, supplemented with the gradient information. Alas, even for a
moderate accuracy of the output approximate eigenvalue, the number of
iterations (and $\sigma_{\min}$ evaluations) of {\tt fminu}, usually
exceeded the number of $100$.

The main part of this paper aims to determine a strategy for selecting the step size.
Using mainly results from our work for the efficient computation of the pseudospectrum
\cite{GK.99}, we give a simple and computationally efficient formula for the step size,
which as we shall see, when incorporated in {\tt sdeig} gives a powerful iteration.

\section{Determining Step Size}\label{sec:STEP.SIZE}

In this section we will use theorems proven in \cite{GK.99}. Denote by
$D^{\circ}(z,r)$ the open disk with center $z$ and radious $r$. We usually
denote with $\sigma_{\min}(z)$ the minimum singular value of $zI-A$ and with
$u_{\min}(z),v_{\min}(z)$ the corresponding left and right singular vectors.

\begin{theorem} \label{th:SMALL.EXCLUSION}
If $\sigma_{\min}(zI-A)=r>\epsilon$, then $D^{\circ}(z,r-\epsilon)\cap
\Lambda_{e}(A) = \emptyset$. In the extreme case $\epsilon=0$, we have $
D^{\circ}(z,r)\cap \Lambda(A) = \emptyset$.
\end{theorem}

\newtheorem{corollary}{Corrolary}

Is $\sigma_{\min}(zI-A)$ a good candidate for step size? The answer is
affirmative in the case of normal matrices, as the following corollary
shows:
\begin{corollary} \label{th:normal.excl} If $N \in \mathbb{C}^{n\times n}$ is
normal, $z_i, i = 1, \ldots, n$ are its eigenvalues, and
$\sigma_{\min}(zI-N)=r>\epsilon$ then the perimeter $\partial
D(z,r-\epsilon)$ contains at least one point of $\Lambda_{\epsilon}(N)$. In
the extreme case $\epsilon=0$, the disk $D(z,r)$ "touches" at least one
eigenvalue of $N$.
\end{corollary}

Corollary \ref{th:normal.excl} along with the observation that pseudospectra
of normal matrices are concentric circles around eigenvalues, guarantee that
{\tt sdeig} always finds in one step  the eigenvalue closest to the starting
point. However, in non-normal matrices the pseudospectra expansion can be
very fast. In other words, the value of $\sigma_{\min}(z)$ can be very small
comparing to the distance of $z$ from $\Lambda(A)$, so that the convergence
of {\tt sdeig} is very slow. Our basic effort will be to measure and cancel
this fast pseudospectra expansion. To this aim we define the
\textbf{pseudospectrum sensitivity} at a point $z$ to be

\newtheorem{definition}{Definition}

\begin{definition}
\begin{eqnarray}
pss(z)=|v_{\min}^{H}(z) *u_{\min}(z)|^{-1}
\end{eqnarray}
\end{definition}

Using the above definition the effectiveness of Theorem
\ref{th:SMALL.EXCLUSION} can be extended as follows.

\begin{theorem} \label{th:LARGE.EXCLUSION}
Let $z=x+iy \in \mathbb{C}\setminus \Lambda(A)$ and $L$ the line segment
between $z$ and the closest eigenvalue. If $pss(z)\leq pss(z')$ for every
point $z' \in L$, then $D^{\circ}(z,\sigma_{\min}(z)\cdot pss(z))\cap
\Lambda_{\epsilon}(A) = \emptyset$.
\end{theorem}
\begin{proof}
Let $g(x,y)=\sigma_{\min}((x+yi)I-A)$. By definition $g(x,y)-g(x_{0},y_{0})
= \oint^{(x,y)}_{(x_{0},y_{0})} \nabla g(x,y)\vec{n}_{L}dL$ , where
$\vec{n}_{L}$ is a normalized vector parallel to $L$. Then using Theorem
\ref{th:gradient} we have
\begin{eqnarray*}
|g(x_0,y_0)-g(x_{1},y_{1})| = |\oint^{(x_1,y_1)}_{(x_{0},y_{0})} \nabla
g(x,y)\cdot \vec{n}_{L} \hspace{0.1cm} dL|  \leq
\oint^{(x_1,y_1)}_{(x_{0},y_{0})} |\nabla g(x,y)\cdot \vec{n}_{L}| ~ dL \\ =
\oint^{(x_1,y_1)}_{(x_{0},y_{0})} pss(z)^{-1} ~dL \leq pss(z_0)^{-1}
\oint^{(x_1,y_1)}_{(x_{0},y_{0})} ~dL = pss(z_0)^{-1} \cdot|z_0-z_{1}|.
\end{eqnarray*}
\end{proof}

Of course, the key question concerning \ref{th:LARGE.EXCLUSION} is its
applicability, given that its assumption does not hold in general. Notice
first that even in cases where the assumption does not hold, the statement
of the theorem can be true, since $pss(z)$ is not in general constant in
$L$. In \cite{GK.99} we showed that the above statement is in general true,
with the exception of some areas between eigenvalues. We could briefly
explain this phenomenon by means of some basic theorems and observations.

\begin{theorem} \label{th:PSS.Eigenvalue}
Let $z_{0}$ be a simple eigenvalue of $A$ with $y,x$ the corresponding left
and right eigenvectors. If we define the eigenvalue condition as
$\kappa(z_{0})= \frac {||{y^{H}}|| \cdot ||{x}|| } {|y^{H}x|} $, then
\begin{eqnarray*}
\lim_{z\rightarrow z_{0}}pss_{A}(z)= \kappa(z_{0})^{-1}
\end{eqnarray*}
\end{theorem}

Observe now that in the limit $|z|\rightarrow +\infty$, we have $v^{H}_n
u_n\rightarrow 1$. This fact along with Theorem \ref{th:PSS.Eigenvalue} and
continuity of $pss(z)$ is an indication of the $pss(z)$ behavior with
$dist(z,\Lambda(A))$ increasing. Of course, our claim gives no guarantee for
the $pss(z)$ behavior close to $\Lambda(A)$. The following theorem taken
from \cite{Stewart.Sun} gives a local description around simple eigenvalues.

\begin{theorem} \label{th:secorder}
Let $z_e$ be a simple eigenvalue of $A$ with $y,x$ the corresponding left
and right eigenvectors and $\tilde{A}=A+E$ be a perturbation of $A$. Then
there is a unique eigenvalue $\tilde{z_e}$ of $\tilde{A}$ such that
\begin{eqnarray}
|\tilde{z_e}-z_e|\leq \frac {||y^{H}|| \cdot ||{x}|| \cdot ||E||}
{|y^{H}x|}+ O(||E^{2}||)
\end{eqnarray}
\end{theorem}

\begin{figure}[h]
\begin{center}
\leavevmode \epsfxsize = 2.0in \epsfysize =2.0in \epsfbox{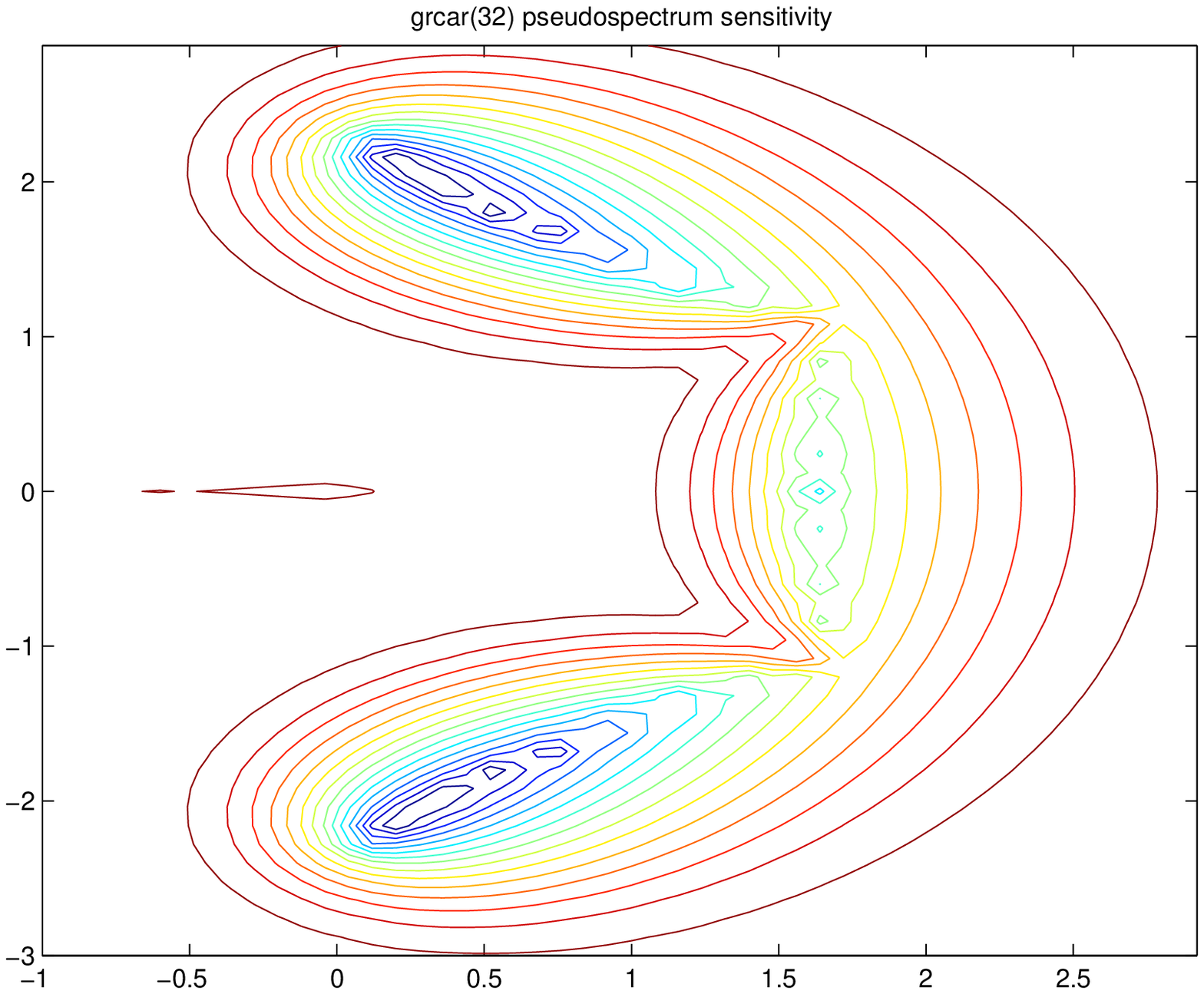} \epsfxsize
=2.0in \epsfysize =2.0in \epsfbox{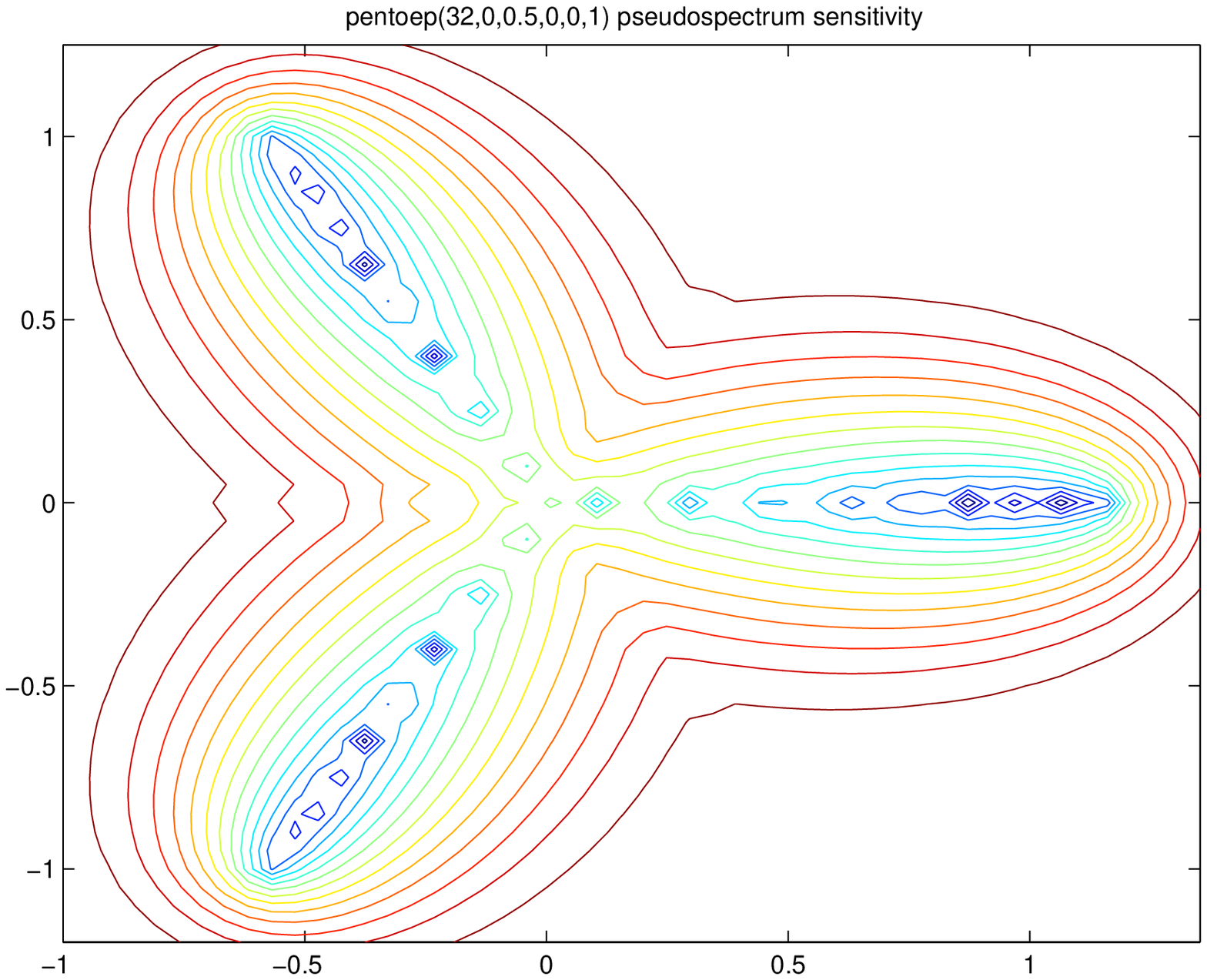}
\end{center}
\caption[]{Contour plots of $pss(z)$ for matrices grcar(32) and
propeller(32) .} \label{fig:pss.contour}
\end{figure}

Second order terms in Theorem \ref{th:secorder} translate directly in a
decrease of $pss(z)$ locally around simple eigenvalues. Concerning areas
between eigenvalues where the statement of Theorem \ref{th:LARGE.EXCLUSION}
fails, we can use a basic result from \cite{LE.98}.
\begin{eqnarray} \label{th:AVOID}
v_{\min}^{H}(z)*u_{\min}(z)=0\Leftrightarrow \exists E \textnormal{ such
that $z$ is a double eigenvalue of $A+E$}.
\end{eqnarray}
In that case $\sigma_n=\sigma_{n-1}$ is a necessary condition. Since
$pss(z)$ is continuous there are areas around $z$ where $pss(z)$ is very
small and thus gives very large incorrect exclusions. We will defer further
discussion on the behavior of pseudospectrum sensitivity behavior until
section \ref{sec:pss.behavior}. The final argument of this section is given
in Figure \ref{fig:pss.contour}, which essentially shows that the assumption
of Theorem \ref{th:LARGE.EXCLUSION} is very weak.

We have set grounds which justify our selection for the step size:
\begin{eqnarray}
step\_size=\sigma_{\min}(z) \cdot pss(z)
\end{eqnarray}

\section{Numerical Results}\label{sec:RESULTS}

In this section we present numerical results that in principal show the
effectiveness of {\tt sdeig}. Starting points have been selected using the
{\tt enclose} algorithm of \cite{GK.99}. We use the complete {\tt svd}
decomposition, but essentially the same results are obtained with Arnoldi
methods, for example with the function {\tt svds} of {\tt Matlab}. In the
following, with $\lambda$ we denote an eigenvalue computed by {\tt eig},
with "eig\_sens" we denote the sensitivity of the corresponding eigenvalue
computed by {\tt eigsens}. All test matrices and function {\tt eigsens} are
from the {\tt Testmatrix Toolbox} for {\tt Matlab}.

\begin{figure}[h]
\begin{center}
\leavevmode \epsfxsize = 2.0in \epsfysize =2.0in \epsfbox{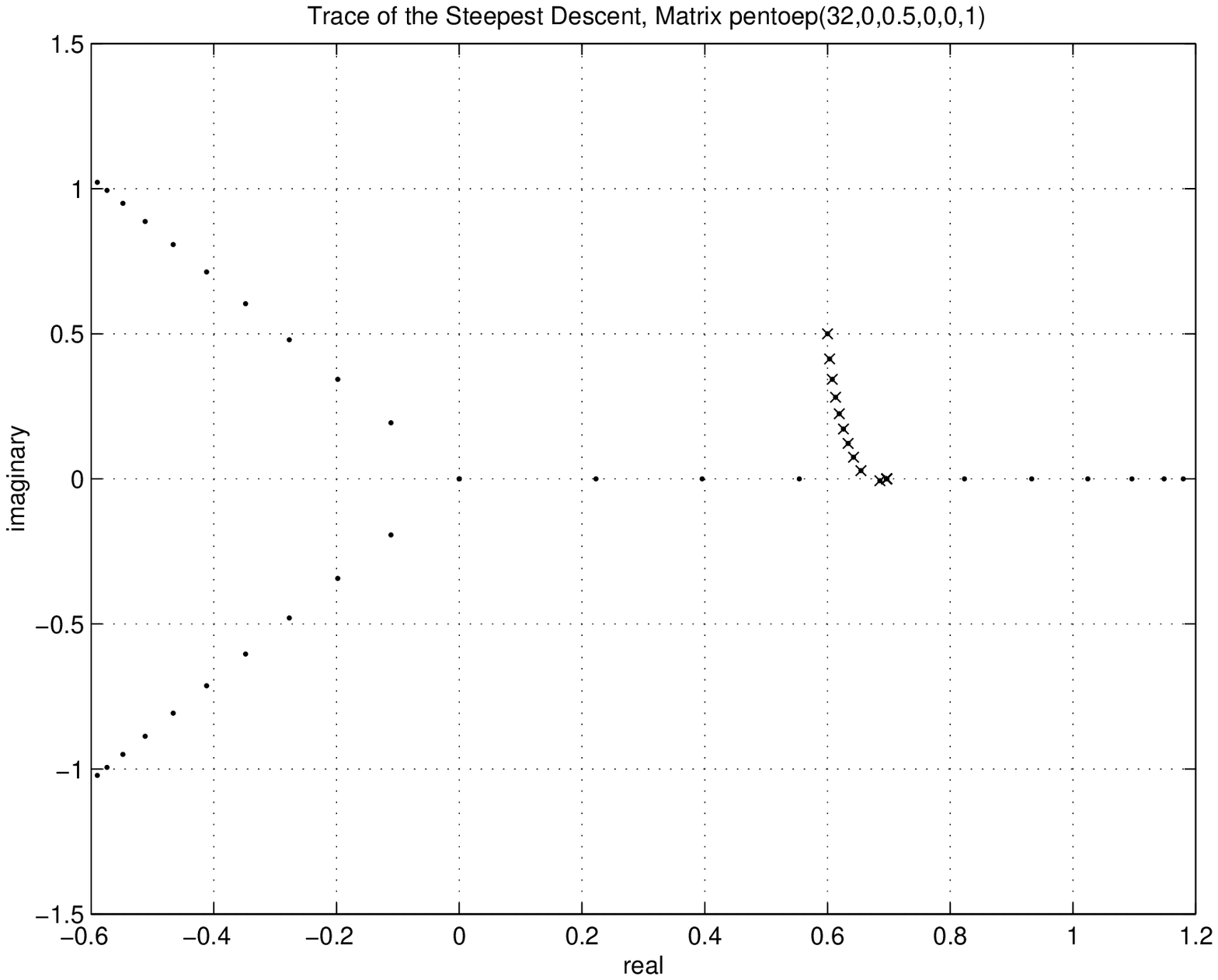} \epsfxsize
=2.0in \epsfysize =2.0in \epsfbox{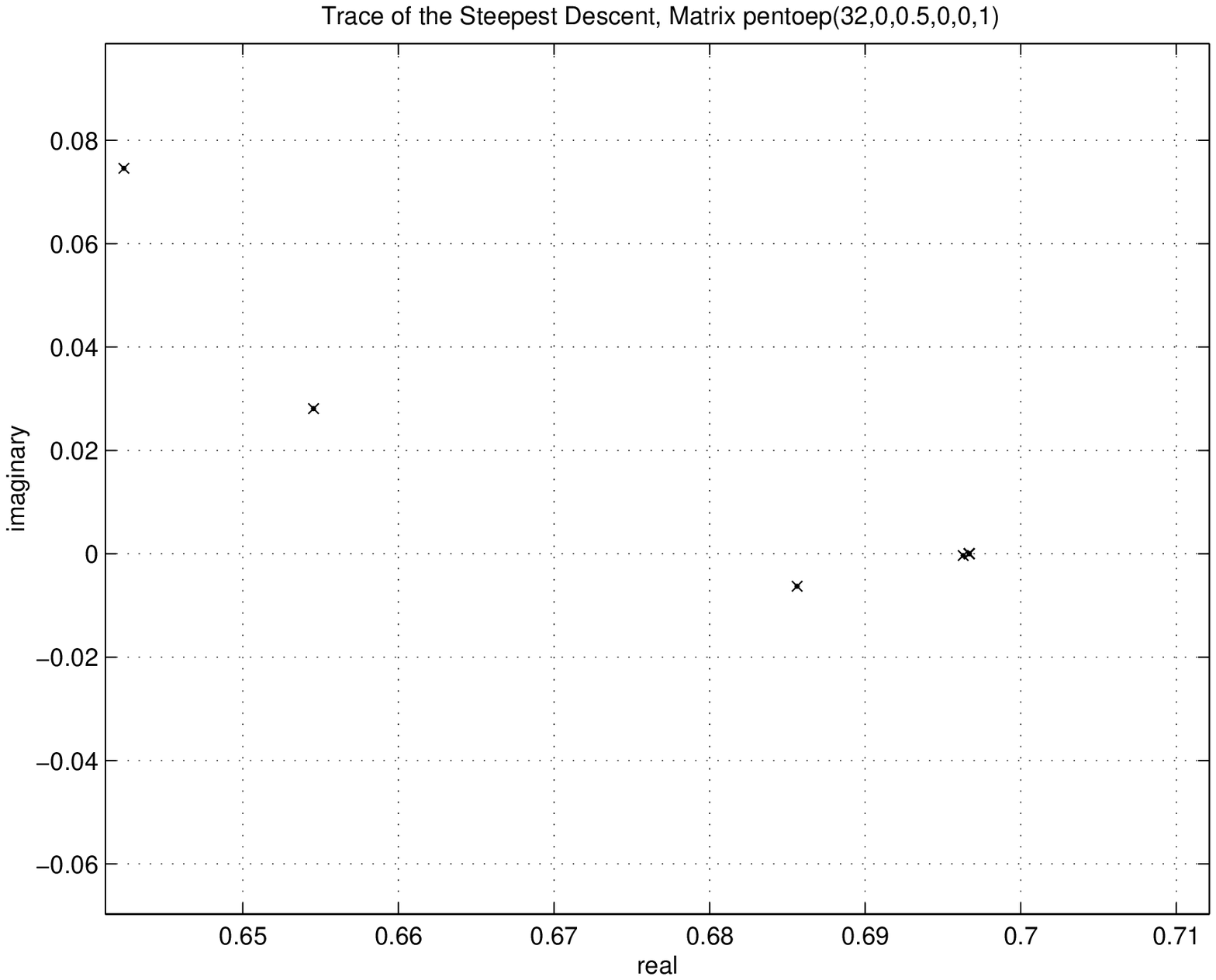}
\end{center}
\caption[]{The steepest descent process}
\label{fig:show.descent}
\end{figure}

The steepest descent is stopped when the first increase of the numerical
values of $\sigma_{\min}$ occurs. In that sense, the implementation of a
stopping criterion is straightforward. However as we can see in Table
\ref{table:prop1} convergence may have been obtained earlier. Hence, a
better implementation of a stopping criterion could stop avoid some of the
last iterations. An example of the steepest descent process and a zoom in
the area of the eigenvalue is given in Figure \ref{fig:show.descent}.

\begin{table}[p]
\begin{center}
\begin{tabular}{|l|l|l|l|l|l|} \hline
  \bf{Starting Point}& \bf{Iterations} & $\lambda$ &$|\hat{\lambda}-\lambda|$
  & s\_min & \bf{eig\_sens}  \\ \hline
  0.6+0.5i & 15 & 6.96686e-01 & 3.84137e-014 & 4.933e-018 & 9.45e+003 \\
  0.6+0.5i & 14 &             & 9.23759e-014 & 1.383e-017 &         \\
  0.6+0.5i & 12 &             & 1.29933e-006 & 1.374e-010 &         \\
 \hline
\end{tabular}
\caption{Matrix pentoep(32,0,0.5,0,0,1)} \label{table:prop1}
\end{center}
\end{table}

\begin{table}[p]
\begin{center}
\begin{tabular}{|l|l|l|l|l|l|} \hline
  \bf{Starting Point}& \bf{Iterations} & $\lambda$ &$|\hat{\lambda}-\lambda|$
  & s\_min & \bf{eig\_sens}  \\ \hline
   1.0+0.1i & 8  & 1.02409e+00 & 1.90412e-014 & 8.191e-019 & 2.75e+004 \\
   1.0+0.5i & 15 & 9.32779e-01 & 1.85407e-014 & 2.091e-019 & 2.32e+004 \\
   1.3+0.0i & 17 & 1.17998e+00 & 5.10702e-015 & 6.897e-019 & 5.21e+003 \\
   1.3+0.0i & 12 &             & 1.35447e-014 & 1.770e-018 &           \\
  -0.4+0.5i & 9  & -(3.48-6.03i)e-01 & 3.57989e-014 & 2.018e-018 & 9.45e+003\\
 \hline
\end{tabular}
\caption{Matrix pentoep(32,0,0.5,0,0,1)} \label{table:pent1}
\end{center}
\end{table}

In Table \ref{table:pent1} , we see that the number of iterations can be
decreased with the choice of a "good" starting point (15/8), and of an
appropriate stopping criterion {\tt tol} (17/12).

\begin{table}[p]
\begin{center}
\begin{tabular}{|l|l|l|l|l|l|} \hline
  \bf{Starting Point}& \bf{Iterations} & $\lambda$ &$|\hat{\lambda}-\lambda|$
  & s\_min & \bf{eig\_sens}  \\ \hline
   1.0+0.1i & 7  & 1.03032e+00 & 4.04639e-022 & 1.586e-018 & 2.82e+002 \\
   1.0+0.5i & 12 &             & 2.22044e-016 & 4.832e-019 &           \\
  -0.4+0.5i & 7  & -(3.44-6.00i)e-01 & 3.51803e-016 & 3.08e-018 & 1.86e+002\\
 \hline
\end{tabular}
\caption{Matrix pentoep(32,0,0.7,0,0,1)} \label{table:pent2}
\end{center}
\end{table}

In the examples of Table \ref{table:pent2} , we use a slightly different
matrix changing a single parameter in function {\tt pentoep}. The spectrum
of the new matrix is of similar form with that discussed on Table
\ref{table:pent1}, but eigenvalue sensitivities are clearly reduced. We use
the starting points of Table \ref{table:pent1}. The effects are clear in the
size of $|\hat{\lambda}-\lambda|$ and more importantly in the number of
iterations.

\begin{table}[p]
\begin{center}
\begin{tabular}{|l|l|l|l|l|l|} \hline
  \bf{Starting Point}& \bf{Iterations} & $\lambda$ &$|\hat{\lambda}-\lambda|$
  & s\_min & \bf{eig\_sens}  \\ \hline
   1.0+0.2i & 10  & 9.32039e-01 & 8.88178e-015 & 2.727e-017 & 9.09e+001 \\
   0.8+0.2i & 13  & 7.03348e-01 & 1.14330e-012 & 3.041e-017 & 3.98e+003 \\
   0.6+0.2i & 15  & 4.94699e-01 & 8.37324e-012 & 2.023e-017 & 7.38e+004 \\
   0.4+0.2i & 17  & 3.47946e-01 & 3.01335e-011 & 1.727e-017 & 1.08e+006 \\
   0.2+0.2i & 18  & 1.72127e-01 & 1.71558e-011 & 1.634e-019 & 3.26e+007 \\
 \hline
\end{tabular}
\caption{Matrix kahan(32)} \label{table:kahan}
\end{center}
\end{table}

In the examples of Table \ref{table:kahan} we use a nice property of matrix
kahan(32).  Its eigenvalues lie on the real axis in the range
$[0.11,1.0.1]$, and their sensitivities satisfy $e_i<e_j \Leftrightarrow
sens(e_i)>sens(e_j)$. We observe that the increase of eigenvalue (and
pseudospectrum) sensitivity reflects in a monotone way in the number of
iterations and in the size of $|\hat{\lambda}-\lambda|$. However, an
O($10^{6}$) increase of the eigenvalue sensitivity does not even double the
number of iterations. We also observe that independently from eigenvalue
sensitivity, $\sigma_{\min}(z-A)<10^{-16}$. For a variety of similar results
refer to \cite{BK.98}.

\section{The role of clustering}\label{sec:CLUSTERS}

In the previous section we observed a connection between the number of
iterations and eigenvalue sensitivity (or closeness of the original matrix
to defective matrices -- see an argument in \cite{Stewart.Sun}). But, how
clusters of eigenvalues on the complex plane can affect the number of
iterations? We settle the question by means of a numerical example, and an
extreme case theorem.

\begin{figure}[h]
\begin{center}
\leavevmode \epsfxsize = 2.0in \epsfysize =2.0in \epsfbox{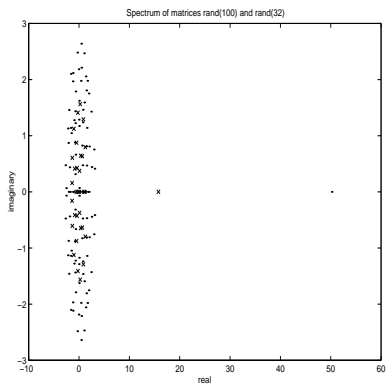}
\end{center}
\caption[]{Examples of geometric clusterings}
\label{fig:clusters}
\end{figure}

In Figure \ref{fig:clusters} we plot spectra of {\tt rand(32)} and {\tt rand(100)}. It is
obvious that eigenvalues of {\tt rand(100)} exhibit a more intense geometric clustering.
However, sensitivity of eigenvalues is less than $100$ in both cases. In particular { \tt
mean(eigsens(rand(32)))=4.59} and { \tt mean(eigsens(rand(100)))=6.25}. For a small
number of starting points inside the cluster, we obtained { \tt
mean(iterations(rand(32)))=6.8} and \\ {\tt mean(iterations(rand(100)))=6.9}. We now give
an extreme case theorem.

\begin{theorem} \label{th:eigeqsvd}
Let $N$ be a normal matrix with eigenvalues
$\lambda_{1},..,\lambda_{n}$, with
$|\lambda_{1}|>..>|\lambda_{n}|$, and singular values
$\sigma_{1}>..>\sigma_{n}$. The following hold:
\begin{enumerate}
\item $|\lambda_{i}|=\sigma_{i}$.
\item For every $z\in \mathbb{C}$, if  $zI-N=U\Sigma V^{H}$,
there is a permutation $\pi$, such that \newline
$\lambda_{\pi(i)}=z-S(i,i)*V(:,i)'*U(:.i)$, $i=1,...,n$.
\end{enumerate}
\end{theorem}
\begin{proof}
For claim 1, see \cite{GJSW.87}. In the case of normal matrices we have
$u^H_iv_i$ for every singular value $\sigma_i$. Let $z_i$ be the eigenvalue
for which $|z_i-z|=\sigma_i(z)$. Let $L$ be the line between $z$ and $z_i$.
For each point $z'$ of $L$, there is some $j<i$ such that
$|z'-z_i|=\sigma_j(z'I-A)$. We now split $L$ in segments where the index $j$
is constant and we consider the integral, as we did in proof of Theorem
\ref{th:LARGE.EXCLUSION} for the corresponding $\sigma_j$ along each
segment. Observe that since $N$ is normal, we have $u_i(z')^Hv_i(z')=1$, for
every $i,z'$, and the inequalities of the proof become strict equalities.
Hence, the gradient of $\sigma_j$ must be parallel to $L$ at every point of
it.
\end{proof}

In Table \ref{table:eig.svd} we give the eigenvalues of a normal
$5 \times 5$ matrix with no trivial eigenvectors, computed using
the algorithm of Theorem \ref{th:eigeqsvd}, for the point $z=1.5$.
We also compare with values computed by {\tt eig}.

\begin{table}[h]
\begin{center}
\begin{tabular}{|l|l|l|l|l|l|} \hline
  {\bf computed eigenvalues}& $|\lambda_{eig}-\lambda|$
  \\ \hline
   4.354989548789193e+000 +2.213138973009386e-001i  & 1.766e-015  \\
   3.902417714058627e+000 +3.619844954770191e+000i  & 3.203e-015  \\
   2.750755557087913e+000 +9.117844570868215e+000i  & 1.432e-014  \\
   3.344861206758000e+000 +9.204041639930084e+000i  & 7.105e-015  \\
   7.136815944761636e+000 +8.562489580579660e+000i  & 1.387e-014  \\
 \hline
\end{tabular} \label{table:eig.svd}
\caption{Eigenvalues by means of \tt{svd}}
\end{center}
\end{table}

A well known result (see \cite{GoVa96}), is that {\tt svd} requires (by a
factor) less flops than {\tt eig}. In our 5x5 example, {\tt flops(eig)=
4695,flops(svd)=2293}. Theorem \ref{th:eigeqsvd} directly suggests a way for
exploiting normality in the computation of eigenvalues of normal matrices
with direct or indirect methods.

This section established also, that the size of the matrix does not affect
the number of iterations.

\section{Using Exclusions}\label{sec:EXCLUIONS}

In the previous sections we established {\tt sdeig}, a steepest
descent algorithm for the computation of a simple eigenvalue. The
central observation of this section is that the steepest descent
procedure can give additional useful information. We repeat the
key claim of Section \ref{sec:STEP.SIZE}.

\newtheorem{claim}{Claim}

\begin{claim}
$D^{\circ}(z,\sigma_{\min}(zI-A)\cdot pss(z))$ and $\Lambda(A)$
have under weak conditions no common points.
\end{claim}

\begin{figure}[h]
\begin{center}
\leavevmode \epsfxsize = 2.5in \epsfysize =2.5in \epsfbox{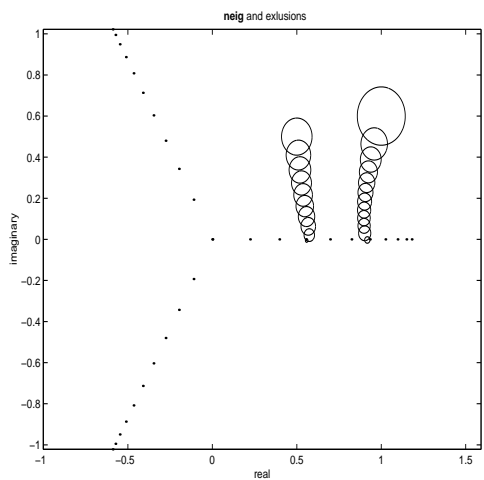}
\end{center}
\caption[]{Using Exclusions}
\label{fig:show.exclusion}
\end{figure}

This claim implies that we can use all the $\sigma_{\min}(z)$
computations of {\tt sdeig} to obtain regions of the complex plan
which do not contain eigenvalues. Figure \ref{fig:show.exclusion}
shows the area excluded after two runs of {\tt eig}. The use of
"localization information" could be an additional feature of an
eigenvalue algorithm which will use our ideas. Of course,
exclusions will be most informative in the case of hermitian and
in general, normal matrices.

In Figure \ref{fig:show.exclusion} also observe that exclusion
disks intersect. A more sophisticated algorithm could choose a
slightly bigger step size (at least in areas away from
eigenvalues) in order to diminish the number of iterations. A
natural way to try such an improvement would be to use second
order information for $\sigma_{\min}$. These are open issues for
future research.

\section{The Oleary-Stewart Result}

The iteration of {\tt sdeig} is given in \cite{OS.98}  slightly disguised.
It is interesting that the iteration is derived from a very different
perspective. The starting point is the following variant of the Rayleigh
quotient method. Let $\lambda$ be a simple eigenvalue of $A$ with right and
left eigenvectors $x$ and $y^H$. Let $\tilde{v}$ and $\tilde{w}^H$ be
approximations to $x$ and $y^H$, and let $\tau$ be an approximation to
$\lambda$. The new approximations $\hat{v},\hat{w}^H$ and $\hat{\tau}$ are
generated as follows:
\begin{eqnarray*}
1. & \hat{v}=(A-\tau I)^{-1} \hat{u}~~~~~
\\ 2. & \hat{w}^H=\hat{w}^H(A-\tau I)^{-1}
\\ 3. & \hat{\tau} = \hat{w}^HA\hat{v}/\hat{w}^H\hat{v}~~~~~~
\end{eqnarray*}
The precedure is then iterated. The basic observation of \cite{OS.98} is
that if $\tau$ is close to an eigenvalue the singular vectors should
approximate the eigenvectors. Using only this observation the Rayleigh
quotient iteration is transformed to the iteration of {\tt sdeig}. Finally,
it is given a proof that the iteration achieves local quadratic convergence.
\begin{theorem} \label{th:OS}
if $\epsilon=|\lambda-\tau|$ is sufficiently small, there is a constant
$\mu$ such that:
\begin{eqnarray}
\hat{\epsilon}\equiv  |\lambda-\hat{\tau}| \leq
C((\epsilon/\mu)^2+(\epsilon/\mu)(||\delta_u||+||\delta_w||)+||\delta_v||||\delta_w||),
\end{eqnarray}
where
\begin{eqnarray}
C=\frac {2k^2||A||}{\sqrt{1-(\epsilon/\mu)^2}
-\kappa(2(\epsilon/\mu)+||\delta_v||+||\delta_w||)}
\end{eqnarray}
The constant $\mu$ is the lower bound of $\sigma_{n-1}(z)$ in a neighborhood
of $\lambda$ and $\delta_u,\delta_w$ denote inaccuracies in the computation
of the singular vectors $u,w$.
\end{theorem}

Unfortunately, the result of \cite{OS.98} holds only for (unbounded) small
neighborhoods around simple eigenvalues, and fails to explain the general
convergence exhibited by {\tt sdeig}. It is however useful in proving
rigorously the ultimate convergence of {\tt sdeig}, and establishing a
better stopping criterion.

\section{Pseudospectrum Sensitivity Behavior} \label{sec:pss.behavior}

Figure \ref{fig:pss.contour} shows that pseudospectrum sensitivity exhibits
a behavior more structured than that implied from theorems and observations
in Section \ref{sec:STEP.SIZE}. In this section we will attempt to give a
more accurate picture of that structure. We first give a general theorem
taken from \cite{Gallestey.Thesis} (see also \cite{Hille.Phillips}).

\begin{theorem} \label{th:sh}
Let $X$ be a Banach space and $\Omega$ an open subset of $\mathbb{C}$.
Consider a function $f: \hat{\Omega}\rightarrow X$ holomorphic in $\Omega$
and continuous in $\tilde{\Omega}$. Then $||f||$ is subharmonic in $\Omega$.
\end{theorem}

We will need one more basic result contained in \cite{Kato.80} (see page
84).

\begin{theorem} \label{th:Kato}
If $A(z) \in \mathbb{C}^{n\times n}$ is an holomorphic function of $z \in
\mathbb{C}$, and $\lambda(z)$ is a simple eigenvalue of $A(z)$ in some
domain $\Omega \subset \mathbb{C}$ then $\lambda(z)$ is holomorphic in
$\Omega$. Moreover, the corresponding left and right eigenvectors $y^H(z),
x$ are holomorphic in $\Omega$.
\end{theorem}

We are now ready to give the main theorem.



\begin{theorem}
If $A \in \mathbb{C}^{n\times n}$, pss(z) is subharmonic in every domain
$\Omega \in \mathbb{C}$, where $\sigma_{\min}(z)$ is simple.
\end{theorem}

\begin{proof}
Let $zI-A=U(z)^H\Sigma(z) V(z)$. It is well known that the
hermitian matrix $(zI-A)^H(zI-a)$ has eigenvalues
$\sigma_1(z)\geq\ldots\geq\sigma_n(z)$, and corresponding left and
right eigenvectors $u_1(z),\ldots,u_n(z)$ and
$v_1(z),\ldots,v_n(z)$. Clearly $\sigma_{n}(z)=\sigma_{\min}(z)$
is simple if $\sigma_{n-1}(z)\neq \sigma_{n}(z)$. In that case,
from Theorem \ref{th:Kato}, $\sigma_n$ and $u_i(z)^H,v_i(z)$ are
holomorphic in $z$ . Hence $u_i^{H}(z)\cdot v_i(z)$ is
holomorphic. Subharmonicity follows from Theorem \ref{th:sh} and
the definition of $pss(z)$.
\end{proof}

The Maximum Principle is an elementary property of subharmonic functions.
Specifically, the following theorem holds.

\begin{theorem} \label{th:MP}
Let $\Omega$ be a bounded domain in $\mathbb{C}$ and $\partial \Omega$ be
its boundary. Suppose $f$ is subharmonic in $\Omega$. Then
\begin{eqnarray}
\sup_{s\in\Omega}f(s) = \max_{\zeta\in\partial\Omega}f(z)
\end{eqnarray}
\end{theorem}

The Maximum Principle can be directly applied to $pss(z)$. Notice also, that
since $pss(z)=0 \Rightarrow \sigma_n=\sigma_{n-1}$, the function
$pss(z)^{-1}$ is subharmonic in the same domain. Hence, the Maximum
Principle can also be applied to $pss(z)^{-1}$. These observations are
summarized in the following corollary.

\begin{corollary} \label{cor:MP.PSS}
Let $\Omega=\mathbb{C}\setminus \{z:\sigma_n(z) \neq \sigma_{n-1}(z)\}$.
Then
\begin{eqnarray}
\nonumber\sup_{s\in\Omega}pss(s) = \max_{\zeta\in\partial\Omega}pss(z) \\
\inf_{s\in\Omega}pss(s) = \min_{\zeta\in\partial\Omega}pss(z)
\end{eqnarray}
\end{corollary}

We can strengthen the result by taking into consideration the behavior of
$pss(z)$ close to points where $\sigma_n(z)=\sigma_{n-1}(z)$. These are
exactly the points where disjoint curves of $\Lambda_{\epsilon}(A)$ collide.

Corollary \ref{cor:MP.PSS} together with some of the observations of Section
\ref{sec:STEP.SIZE} plausibly explains why $pss(z)$ decreases in a monotonic
way as $dist(z,\Lambda(A))$ increases. However, the similarity of $pss(z)$
contours with $\Lambda_{\epsilon}$ contours is striking and obviously more
structured. We will need the extension of Sun's Theorem.

\begin{theorem}
Let $U(x,y)\Sigma V^H(x,y)=A-(x+iy)I$, where $\sigma=\sigma_n(x_0,y_0)$ is
simple and non-zero. Let $\tilde{U}=(u_1 u_2 \ldots u_{n-1})$,
$\tilde{V}=(v_1 v_2 \ldots v_{n-1})$ and $\tilde{\Sigma} =
diag(\sigma_1,\ldots,\sigma_{n-1})$. We then have
\begin{eqnarray*}
\frac{\partial^2 \sigma}{\partial x^2}(x_0,y_0) &=& Re\{r^H\Phi r + l^H\Phi
l +2l^H \Psi r\}+Im\{u^H_nv_n\}^2/\sigma_n,
\\ \frac{\partial^2\sigma}{\partial x \partial y}(x_0,y_0) &=& -Im\{2l^H\Psi r\}+
Im\{u_n^Hv_n\}Re\{u_n^Hv_n\}/\sigma_n,
\\ \frac{\partial^2 \sigma}{\partial y^2}(x_0,y_0) &=& Re\{r^H\Phi r + l^H\Phi
l -2l^H \Psi r\}+Re\{u^H_nv_n\}^2/\sigma_n,
\end{eqnarray*}
where $\Phi = \sigma_n(\sigma_n^2I-\tilde{\Sigma}^2)^{-1}, \Psi =
\tilde{\Sigma}(\sigma_n^2I-\hat{\Sigma}^2)^{-1}, r=\tilde{U}^Hv_n$ and
$l=\tilde{V}^Hu_n$.
\end{theorem}

\begin{figure}[h]
\begin{center}
\leavevmode \epsfxsize = 2.0in \epsfysize =2.0in \epsfbox{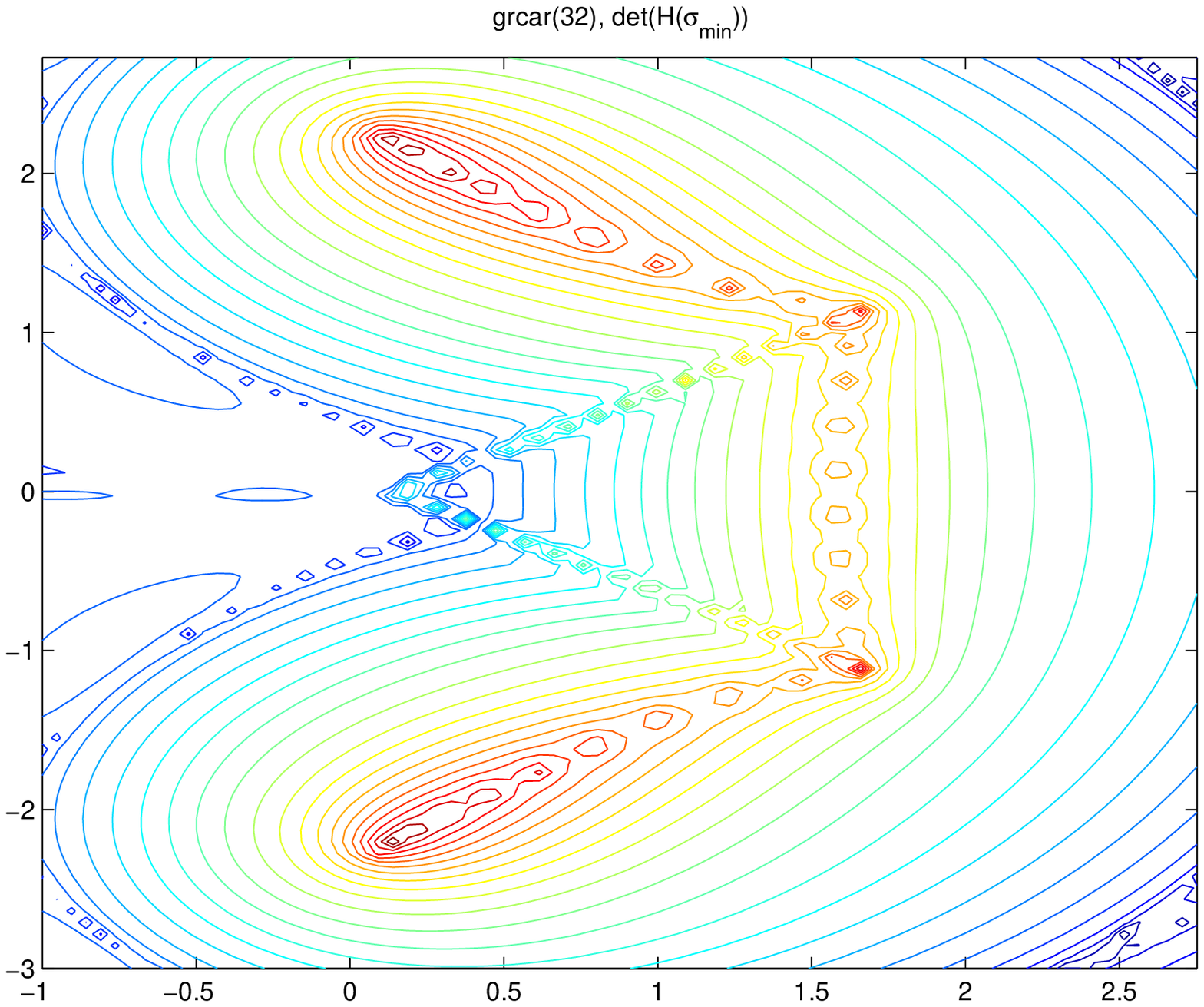} \epsfxsize
=2.0in \epsfysize =2.0in \epsfbox{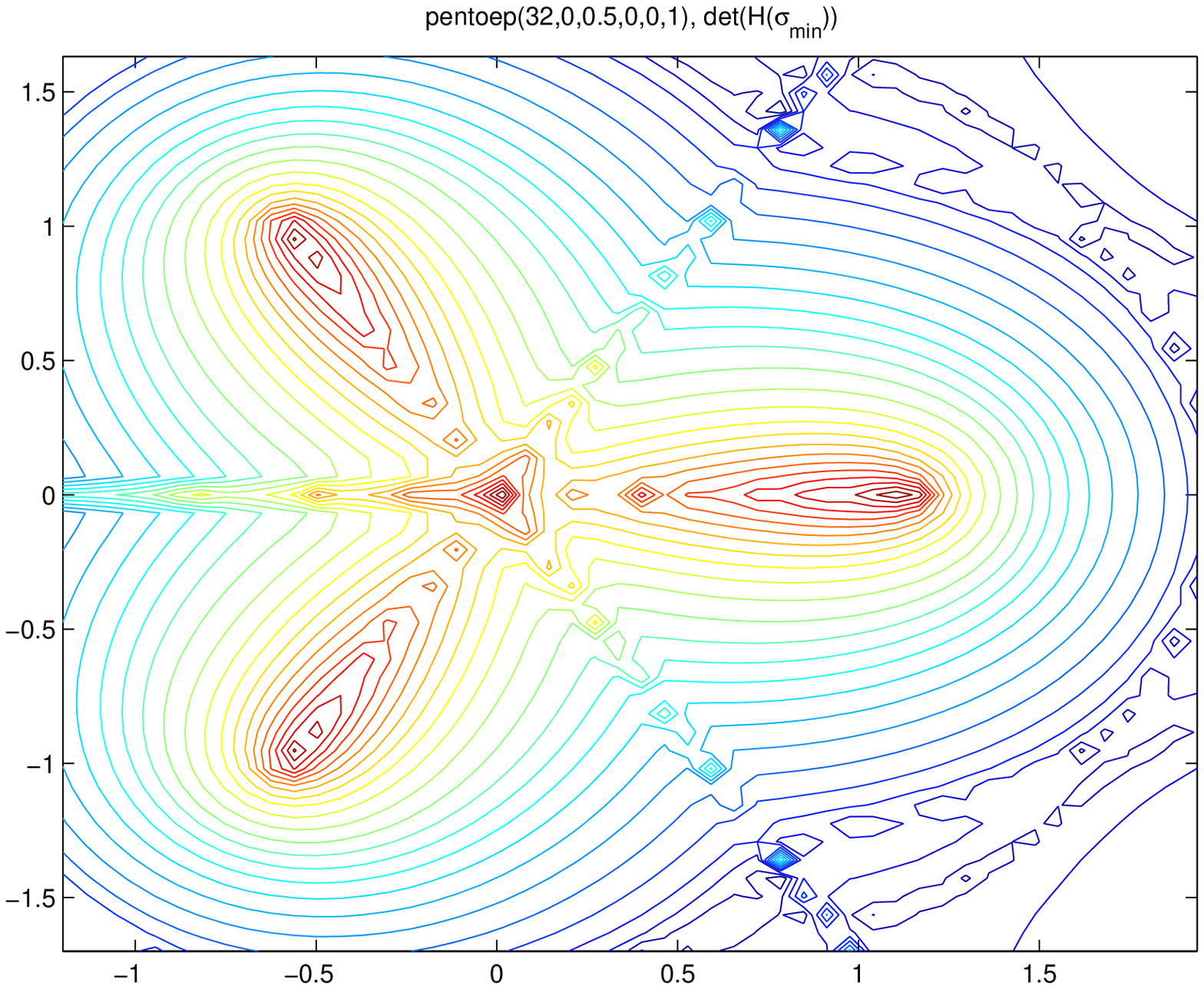}
\end{center}
\caption[]{Contour plots of $det(H(\sigma_{\min}(z)))$ for matrices
grcar(32) and propeller(32) .} \label{fig:hessian}
\end{figure}

Unlike the first order case, computation of second order derivatives is
prohibitive for practical purposes. However, it gives a very good tool for
further theoretical investigation. In Figure \ref{fig:hessian} we plot the
determinant of the Hessian of $\sigma_{\min}(z)$.

Obviously, the pattern is similar with that of $\Lambda_{\epsilon}$ and
$pss(z)$, but clearly, we cannot use a subharmonicity argument. However, we
can observe that expansion is slower than that of $pss(z)$. Intuitively, we
could expect that the sensitivity of analogous constructions for derivatives
of order $k$ asymptotically approaches 1 when $k\rightarrow \infty$. A
quantity that can naturally capture details about these expansion phenomena
is the {\bf pseudospectrum area}, defined in the obvious way.

\begin{figure}[h]
\begin{center}
\leavevmode \epsfxsize = 2.0in \epsfysize =2.0in \epsfbox{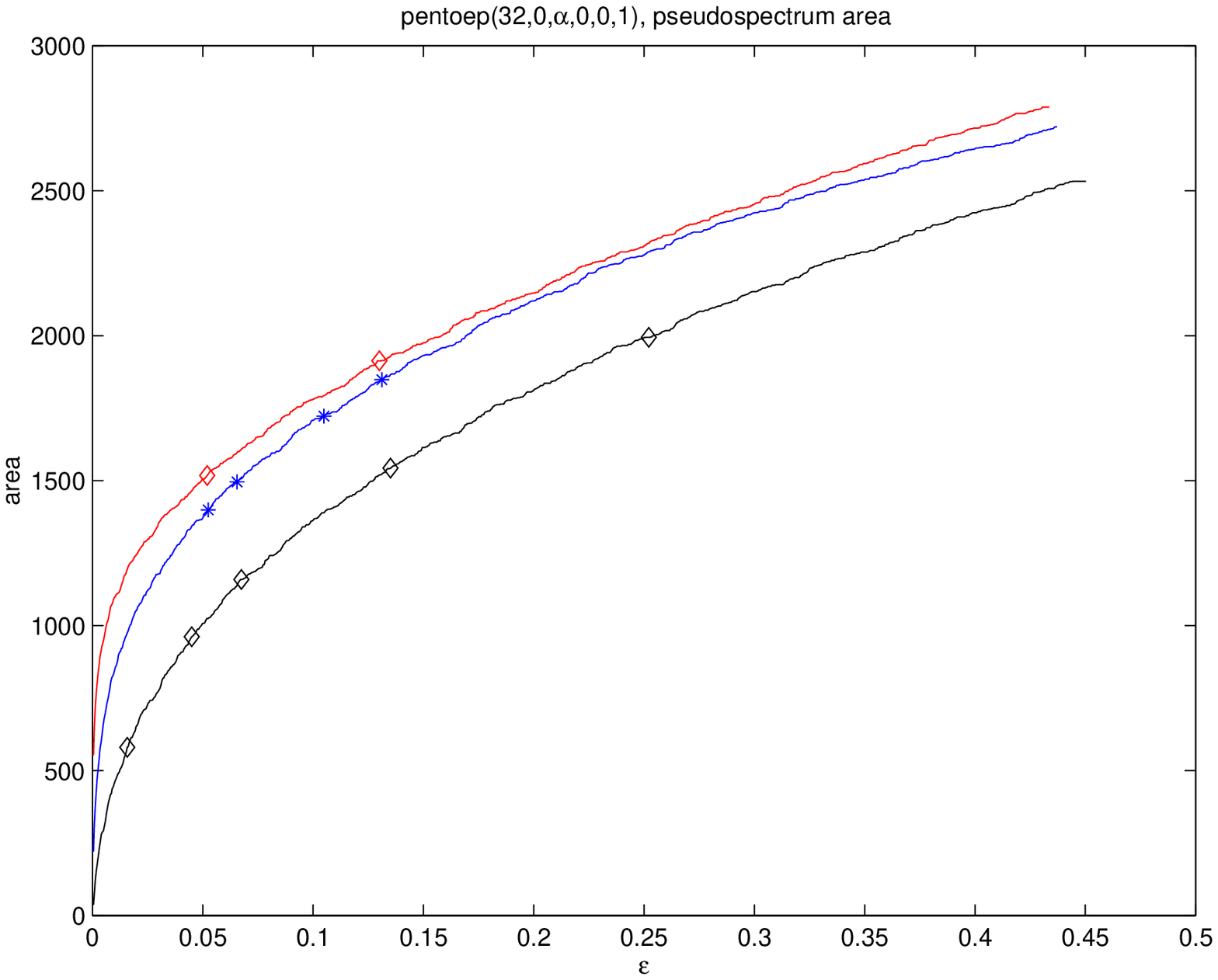} \epsfxsize = 2.0in
\epsfysize =2.0in \epsfbox{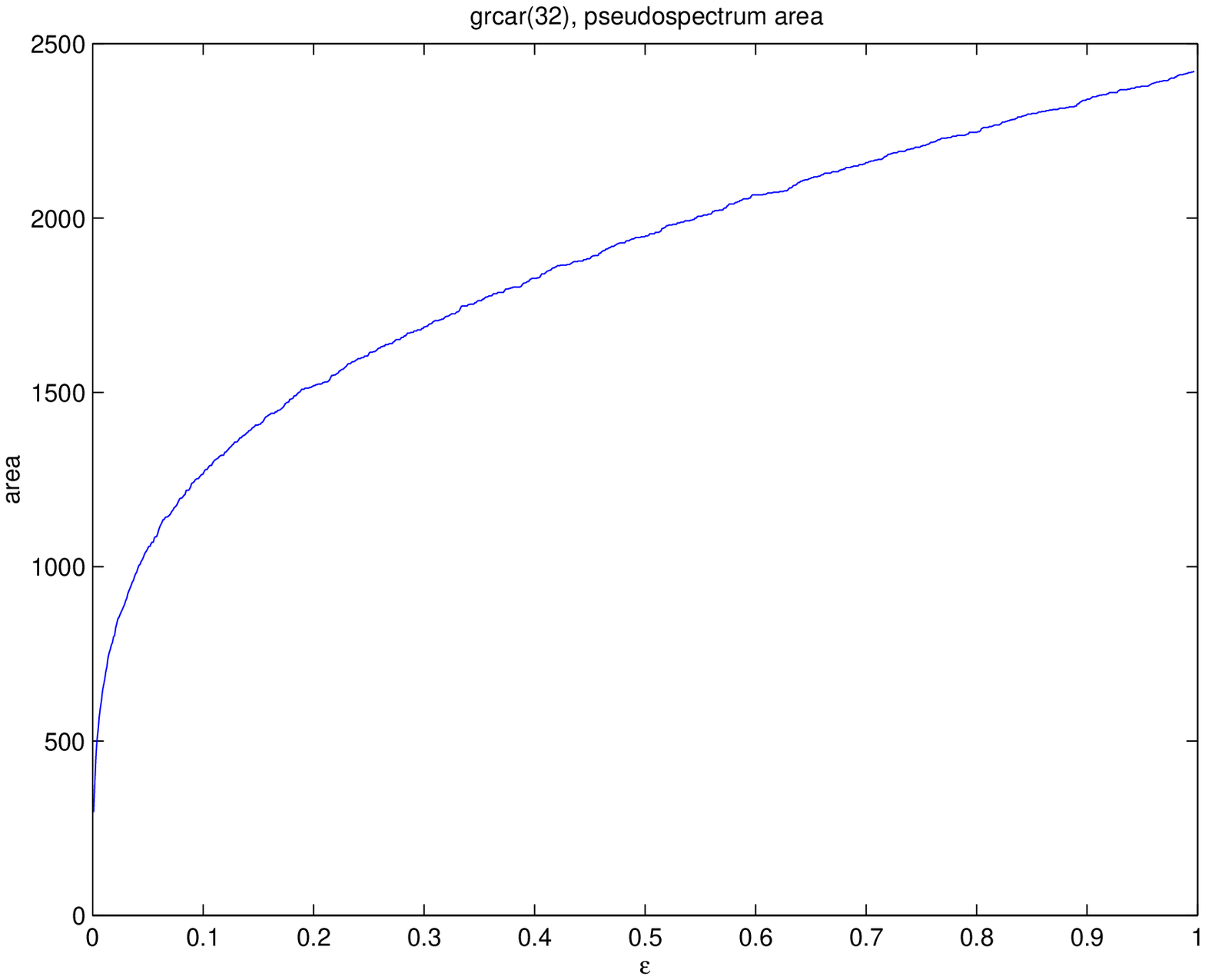}
\end{center}
\caption[]{Evolution of pseudospectrum area} \label{fig:area}
\end{figure}

In Figure \ref{fig:area} we measure the area (as the number of grid points)
of $\Lambda_{\epsilon}$ as a function of $\epsilon$, of matrix $\mathtt{
pentoep}(32,0,\alpha,0,0,1)$, for three different values of parameter
$\alpha$. This family of matrices has the nice property that the spectra are
almost identical for different values of $\alpha$, but the sensitivities of
eigenvalues change significantly. Thus, the $\Lambda_{\epsilon}$ area
depends mainly on $pss(z)$. We also measure the area of $\Lambda_{\epsilon}$
for the matrix ${\mathtt grcar(32)}$. It is clear that a general "inverse
exponential law" shows up, in the sense that despite the initial explosion,
the expansion of pseudospectra stabilizes quickly. The curves are slightly
more steep for bigger sensitivities of eigenvalues.

We now give some theorems which give a way for reasoning about the
pseudospectra expansion and its relation to $pss(z)$. First, a general
theorem of Konrod.

\begin{theorem} \label{th:Konrod}
Let $G\subset\mathbb{R}^{n}$ be a domain and $f:G\rightarrow\mathbb{R}$ a
smooth function, $E_t=\{x\in G|f(x)=t\}$ the level set of the function f and
$ds$ the $(n-1)$-dimensional surface element on $E_t$. Then
\begin{eqnarray*}
meas(G) =\int_{\min f}^{\max f}(\int_{E_{t}}\frac{ds}{|\nabla f|})dt
\end{eqnarray*}
\end{theorem}

Theorem \ref{th:Konrod} directly applies to our case giving the following
corollary.

\begin{corollary} Let $\partial \Lambda_{t}$ be the boundary of
$\Lambda_{t}$ and $dz$ be the 1-dimensional surface element on
$\partial\Lambda_{t}$. Then
\begin{eqnarray*}
area( \Lambda_{\epsilon}(A)) =\int_{0}^{ \epsilon}(\int_{\partial
\Lambda_{t} }(pss(z)dz))dt
\end{eqnarray*}
\end{corollary}

\newpage
\bibliography{../../../Bibliography/J-koutis}

\end{document}